\documentclass{article}
\usepackage{amssymb}
\usepackage{amsmath}
\usepackage[english]{babel}
\usepackage{graphicx}
\usepackage{theorem}

\newtheorem{theorem}{Theorem}

\newtheorem{lemma}[theorem]{Lemma}

\newtheorem{remark}[theorem]{Remark}
\begin{document}
\title{Higher dimensional nonclassical eigenvalue asymptotics}
\author{Brice Camus.\\ Ludwig Maximilians Universit\"at
M\"unchen,\\
Mathematisches Institut, Theresienstr. 39 D-80803
M\"unchen.\\Email: brice.camus@uni-due.de\medskip\\
Nils Rautenberg.\\Ruhr-Universit\"at Bochum, Fakult\"at f\"ur
Mathematik,\\ Universit\"atsstr. 150, D-44780 Bochum, Germany.\\
Email : Nils.Rautenberg@ruhr-uni-bochum.de\\
}
\date{Revised: \today}
\maketitle

\begin{abstract}
\noindent  In this article we extend B. Simon's construction and
results \cite{Sim2} for leading order eigenvalue asymptotics
to $n$-dimensional Schr\"odinger operators with non-confining potentials given by:
$H^\alpha_n=-\Delta +\prod\limits_{i=1}^n |x_i|^{\alpha_i}$ on
$\mathbb{R}^n$ ($n>2$), $\alpha:=(\alpha_1,\cdots,\alpha_n)\in
(\mathbb{R}_{+}^*)^n$. We apply the results to also derive the
leading order spectral asymptotics in the case of the Dirichlet
Laplacian $-\Delta^D$ on domains $\Omega^\alpha_n=\{x\in\mathbb{R}^n: \prod\limits_{j=1}^n |x_j|^{\frac{\alpha_j}{\alpha_n}}<1 \}$.\medskip
\end{abstract}
\section{Introduction and main results.}
Since the seminal work of Weyl \cite{Weyl} and its
generalizations, the eigenvalue asymptotics of the Laplacian
$-\Delta$ on compact domains $\Omega \subset \mathbb{R}^n$ with
various boundary conditions have been understood to encode
information about the geometry of the domain. Let $(\lambda_j)_{j
\in \mathbb{N}}$ denote the sequence of eigenvalues of $-\Delta$
on such a domain endowed with Dirichlet boundary conditions. Let
furthermore $N(E)=\#\{ j : \lambda_j\leq E \}$ be the counting
function of eigenvalues. Not even any regularity of the boundary
$\partial \Omega$ is required for the Weyl law (see \cite{BS}):
\begin{equation*}
N(E)=\frac{\mathrm{vol}(S^{n-1})}{(2\pi)^n}\mathrm{vol}(\Omega)E^{n/2}+\mathrm{o}(E^{n/2}),
\quad \mathrm{as} ~E \to \infty.
\end{equation*}
\noindent This, and further research on additional terms in the
asymptotic expansion
led to the famous question of M. Kac: "Can one hear the shape of a drum?" \cite{Kac}.\\
However, despite the appealing simplicity of this leading order
asymptotics, neither compactness nor finite volume of $\Omega$ are
necessary conditions for purely discrete spectrum of the Dirichlet
Laplacian, denoted from now on by $-\Delta^D$. A class of
$2$-dimensional examples for infinite volume domains with discrete
spectrum, as well as their leading order eigenvalue asymptotics
was given by B. Simon in \cite{Sim2}. He considered 2-dimensional
domains of the form:
\begin{equation*}
\Omega_\alpha:=\{(x,y)\in\mathbb{R}^2 : |x|^\alpha|y|< 1\}
\end{equation*}
with $\alpha>0$ and derived the asymptotics:
\begin{equation*}
N(E)=\begin{cases} \zeta(\alpha)
(\frac{\pi}{2})^{-\alpha}\frac{\Gamma(\frac{1}{2}\alpha+1)}
{\sqrt{\pi}\Gamma(\frac{1}{2}\alpha+\frac{3}{2})}E^{\frac{\alpha+1}{2}}
+ \mathrm{o}(E^\frac{\alpha+1}{2}), \quad \alpha>1\\
\frac{1}{\pi}E\ln{E}+\mathrm{o}(E\ln(E)), \quad \alpha=1,
\end{cases}
\end{equation*}
\noindent where $\zeta$ denotes the Riemannian zeta function and $\Gamma$ denots the Gamma function. For $1>\alpha>0$ the first formula holds if one replaces
$\alpha$ by $\alpha^{-1}$. This article is concerned with an
extension of this example of eigenvalue asymptotics to higher
dimensions. To this end we will determine the spectral asymptotics
for (most members of) the class of Schr\"odinger operators given
by:
\begin{equation*}
H^\alpha_n=-\Delta +\prod\limits_{i=1}^n
|x_i|^{\alpha_i}=-\Delta+V(x), \quad (x_1,...,x_n) \in
\mathbb{R}^n, \quad \alpha_i \in \mathbb{R}_+^* .
\end{equation*}
\noindent If $n \neq 1$ their potentials are non-confining, yet
their spectrum is purely discrete as we will see shortly. In turn
this will enable us to find the eigenvalue asymptotics of the
Dirichlet Laplacian $-\Delta^D$ on domains of the form:
\begin{equation*}
\Omega^\alpha_n=\{x\in\mathbb{R}^n: \prod\limits_{j=1}^n
|x_j|^{\alpha_j/\alpha_n}<1 \}.
\end{equation*}

\paragraph{Remark:}Define the classical energy surfaces:
\begin{equation*}
\Sigma_E=\{ (x,\xi)\in T^* \mathbb{R}^n : ||\xi||^2+V(x) =E\},
\end{equation*}
equipped with the (classical flow invariant) Liouville measures
$\mathrm{dLVol}_E$, and the classical areas:
\begin{equation*}
A(E)=\bigcup\limits_{e\leq E} \Sigma_e=\{ (x,\xi)\in
T^*\mathbb{R}^n: ||\xi||^2 +V(x)\leq e\},
\end{equation*}
for the usual Lebesgue measure $d\mu$. Then, for any choice of the
powers $\alpha_i$, we deal with \textbf{non-compact energy
surfaces of infinite volume}, i.e.
$\mathrm{Lvol}_E(\Sigma_E)=+\infty$. Since we also have
$\mu(A(E))=+\infty$ for $E>0$, Weyl's law does apply, neither in
the classical nor the micro-local category.

\noindent Despite this fact, B. Simon determined the asymptotics
of the Schr\"odinger operators above in the case $n=2$
\cite{Sim2}. We will extend the calculation of these asymptotics
to dimensions $n>2$ covering both the generic case $\alpha_i\neq
\alpha_j, ~\forall i,j$ as well as the most singular case
$\alpha_i=\alpha_j ~ \forall i,j$. Up to a permutation of
coordinates we can freely assume $\alpha_1 \geq \cdots \geq
\alpha_n>0$. Our main results then read as follows:
\begin{theorem} \label{main result}
For $n\geq 2$, assume $\alpha_1 > \cdots > \alpha_n>0$ and define:
\begin{equation*}
d_n=d_n(\alpha):=\frac{\alpha_1+...+\alpha_{n-1}+2}{2\alpha_n}.
\end{equation*}
Then, as $t$ tends to $0^+$, we have:
\begin{equation*}
\lim_{t\rightarrow 0^+} \, t^{(d_n+1/2)} \,
\mathrm{Tr}(e^{-tH^\alpha_n})= \frac{\mathrm{Tr}((H_{n-1}^\alpha)^{-d_n})}{\pi^{n/2}}
\Gamma(d_n+1),
\end{equation*}
where:
\begin{equation*}
\mathrm{Tr}((H^\alpha_{n-1})^{-d_n})=\mathrm{Tr}((-\Delta_{x_1...x_{n-1}}+\prod\limits_{i=1}^{n-1}|x_i|^{\alpha_i})^{-d_n})<\infty,
\end{equation*}
is the spectral-zeta function, evaluated at $d_n$, of the
$(n-1)$-dimensional Schr\"odinger operator obtained by removing
the direction of smallest decay at infinity.
\end{theorem}
\noindent Assuming that all exponents are equal we prove:
\begin{theorem} \label{main result2}
If $\alpha_{1}=\cdots =\alpha_n=\alpha_0$, and thus $d_n=\frac{n-1}{2}+\alpha_0^{-1}$, then
as $t$ tends to $0^+$ we have:
\begin{equation*}
\lim_{t\rightarrow 0^+} \, t^{n/2+\alpha_0^{-1}}|\ln(t)|^{-(n-1)}\,
\mathrm{Tr}(e^{-tH^{\alpha}_n})= \frac{\Gamma(1+\alpha_0^{-1})(n/2+\alpha_0^{-1})^{(n-1)}}{\pi^{n/2}(n-1)!}
.
\end{equation*}
\end{theorem}
\noindent Subsequently, using the Tauberian theorem of Karamata,
cf. \cite{Sim2}, we can prove that the eigenvalue counting
functions $N_{H^\alpha_n}$ given by:
\begin{equation*}
N_{H^\alpha_n}(E)=\# \{j\in\mathbb{N}: \lambda_j\in
\sigma(H^\alpha_n), \lambda_j \leq E\}
\end{equation*}
satisfy the following asymptotic laws as $E\rightarrow \infty$:
\begin{theorem} \label{main result3}
For $n\geq 2$, and assuming $\alpha_1 > \cdots > \alpha_n>0$ we have:
\begin{equation*}
\lim_{E\rightarrow \infty} \, E^{-(d_n+1/2)} \,
N_{H^\alpha_n}(E)= \frac{\mathrm{Tr}((H_{n-1}^\alpha)^{-d_n})\Gamma(d_n+1)}{\pi^{n/2}\Gamma(d_n+3/2)},
\end{equation*}
where $d_n$ and $\mathrm{Tr}((H^\alpha_{n-1})^{-d_n})$ are as in Theorem 1.
\end{theorem}
\begin{theorem} \label{main result4}
If $\alpha_{1}=\cdots =\alpha_n=\alpha_0$, then we have:
\begin{equation*}
\lim_{E\rightarrow \infty} \, E^{-(n/2+\alpha_0^{-1})}\ln(E)^{-(n-1)}
N_{H^{\alpha}_n}(E)= \frac{\Gamma(1+\alpha_0^{-1})(n/2+\alpha_0^{-1})^{(n-1)}}{\Gamma(n/2+\alpha_0^{-1}+1)\pi^{n/2}(n-1)!}
.
\end{equation*}
\end{theorem}
\noindent Finally, this will imply the following result for the
spectrum of the Dirichlet Laplacian $-\Delta^D$ on the domains
$\Omega_n^\alpha$:
\begin{theorem}
For $\alpha_1> \cdots
>\alpha_n>0$ the Dirichlet-Laplacian $-\Delta^D=-\Delta^D_{\alpha,n}$
attached to the domain $\Omega^\alpha_n$, $n\geq 2$ has discrete spectrum and the counting function of eigenvalues satisfies:\\
\begin{equation*}
\lim\limits_{E \to \infty}E^{-(q(\alpha)+\frac{1}{2})}N_{-\Delta^D_{\alpha,n}}(E) = \frac{\mathrm{Tr}((-\Delta^D_{\alpha,n-1})^{-q(\alpha)})\Gamma(q(\alpha)+1)}{\pi^{n/2}\Gamma(q(\alpha)+3/2)}  .
\end{equation*}
where $-\Delta^D_{\alpha,n-1}$ is the Dirichlet Laplacian on the $n-1$-dimensional domain obtained by projecting $\Omega^\alpha_n$ onto the $x_n=0$ hyperplane and
\begin{equation*}
q(\alpha)=\frac{\alpha_1+...+\alpha_{n-1}}{2\alpha_n}.
\end{equation*}
\end{theorem}

\begin{theorem}
The Dirichlet-Laplacian $-\Delta^D$ on $\Omega:=\{ x \in \mathbb{R}^n ~|~ \prod\limits_{i=1}^n|x_i|<1 \}$, $n \geq 2$ has discrete spectrum and:
\begin{equation*}
\lim\limits_{E \to \infty} \ln(E)^{-(n-1)} E^{-\frac{n}{2}} N_{-\Delta^D} (E) = \frac{n^{n-1}}{\Gamma(\frac{n}{2})\sqrt{\pi^n}2^{n-1}(n-1)!}.
\end{equation*}
\end{theorem}

\noindent With regards to other work on similar operators, we first mention
the work of D. Robert \cite{Rob} on the family of potentials:
\begin{equation*} V(x,y)=(1+x^2)^ry^{2l}, \quad r,l
\in \mathbb{R^+}.
\end{equation*}
\noindent Using methods of microlocal analysis, he finds the
leading order asymptotics for the corresponding Schr\"odinger
operators. This precedes Simon's work by about two years and shows
up the great difficulties to derive spectral estimates in presence of a non-confining potential.\\

\noindent Some extensions and alterations of the Simon example using the methods developed in his paper have also been published.
In \cite{AN}, Aramaki and Nurmuhammad have considered the potentials:

\begin{equation*}
V(z)=V(x,y)=||x||^{2p}||x||^{2q}, \quad p,q>0,
\end{equation*}
\noindent and derived the leading order asymptotics for these
cases. More recently, in \cite{EB} Exner and Barseghyan showed the
discreteness of the spectrum, and obtained some bounds on moments
of eigenvalues, for Sch\"odinger operators with potentials
unbounded from below:
\begin{equation*}
|xy|^p -\lambda(x^2+y^2)^{\frac{p}{p+2}},\, p\geq 1,\, \lambda\in
]0,\lambda_ {\mathrm{crit}}[.
\end{equation*}    \\

\noindent The main contribution of our work is to increase dimensions with even more
separate factors in the definition of the potentials (resp. domains for the
Dirichlet Laplacian) while keeping track of the top order coefficients of
the asymptotic expansion of the counting function of eigenvalues.\\

\noindent The constant involved in these asymptotics, explicitely given as a
spectral zeta-function of a lower dimensional operator, is of great
interest and connects very much in line with the results of all the other
work mentioned above.\\

\noindent This work is part of the second author's PhD thesis.\\

\noindent Before we begin our analysis we recall that for a
Schr\"odinger operator with a positively homogeneous potential of
degree $p\neq -2$, i.e. $V(tx)=|t|^p V(x)$ for all $x$, we have
the following scaling relations for all $c>0$:
\begin{gather}
\sigma(-\Delta+c V)=c^{\frac{2}{p+2}} \sigma (-\Delta +V),\label{scale potential}\\
\sigma(-c \Delta +V)= c^{\frac{p}{p+2}}\sigma
(-\Delta+V)\label{scale laplace}.
\end{gather}
These are equalities between spectra and these scaling relations
are independent of the dimension $n$.\\
Following the strategy of B. Simon, who approached the
two-dimensional case with what he calls the ' sliced-bread
inequalities ', we recall now for the readers convenience the classical Tauberian
theorem relating the small time behavior of the quantum partition
function
\begin{equation*}
Z_Q(t)= \mathrm{Tr}\, (e^{-tH}),
\end{equation*}
and the large energy behavior of the counting function
$N_{H}(E)$ as seen in \cite{Sim2}:
\begin{theorem} \textbf{(Karamata's Tauberian Theorem.)}\\
Let $H=-\Delta+V$ with $V$ continuous and non-negative. When $l$
and $d$ are positive, we have the equivalences:
\begin{gather*}
\lim\limits_{E \rightarrow +\infty} E^{-l}N_{H}(E)=c
\Leftrightarrow \lim\limits_{t\rightarrow 0^+} t^{l}
\mathrm{Tr}(e^{-tH})=c \Gamma (l+1),\\
\lim\limits_{E \rightarrow +\infty} \frac{E^{-l}}{(\log E)^d}
N_{H}(E)=c \Leftrightarrow \lim\limits_{t\rightarrow 0^+} \frac
{t^{l}}{|\log(t)|^d}  \mathrm{Tr}(e^{-tH})=c \Gamma (l+1).
\end{gather*}
\end{theorem}

\noindent So we will concentrate on the short time asymptotics of
$Z_Q$ but the infinite volume of the energy surfaces provides a
serious obstacle. As we remarked earlier, in the case of the
operators we are interested in the common trick to estimate the
quantum partition function by a classical integral fails for
precisely this reason. Indeed if we define:
\begin{equation*}
Z_{cl}(t)=\frac{1}{(2\pi)^n}\int\limits_{T^*\mathbb{R}^n}
e^{-t(||\xi||^2+V(x))}dxd\xi,
\end{equation*}
the inequality $Z_Q(t)\leq Z_{cl}(t)$ is valid (this can be viewed
as a convexity property of the exponential function or a consequence of the
abstract Golden-Thompson inequality), but for the product potential:
\begin{equation*}
V_\alpha(x):= \prod\limits_{i=1}^n |x_i|^{\alpha_i},
\end{equation*}
we get:
\begin{gather*}
(2\pi)^n Z_{cl}(t)= (\frac{\pi}{t})^{\frac{n}{2}}
\int\limits_{\mathbb{R}^n} e^{-t V_\alpha(x)} dx\\
=(\frac{4\pi}{t})^{\frac{n}{2}}\Gamma(\frac{\alpha_n+1}{\alpha_n})
\int\limits_{\mathbb{R}_+^{n-1}} \left( t
\prod\limits_{j=1}^{n-1}|x_j|^{\alpha_j}\right)^{-\frac{1}{\alpha_n}}dx_1\dots
dx_{n-1}\\
=(\frac{4\pi}{t})^{\frac{n}{2}}\Gamma(\frac{\alpha_n+1}{\alpha_n})
t^{\frac{1-n}{\alpha_n}}
\prod\limits_{j=1}^{n-1}\int\limits_{0}^{\infty} x_j
^{-\frac{\alpha_j}{\alpha_n}} dx_j= +\infty,
\end{gather*}
none of these integrals being convergent on $[0,\infty]$, for any
choice of the $\alpha_j$ (the convergence
at the origin implying the divergence at infinity and vice-versa).
Despite the fact that the classical
estimate does not yield any useful information, we can easily obtain:
\begin{lemma} \label{easy lemma}
For any $n$ and any ${\alpha} \in (\mathbb{R}_+^*)^n$, the spectrum
of $H^{\alpha}_n$ is discrete.
\end{lemma}
Let us give first a rough proof that the $H^\alpha_n$'s have
discrete spectrum. This approach will not even catch the right
asymptotic power. In dimension 2, this was given by B. Simon in
\cite{Sim1} as the most elementary proof between 6 different
proofs of the
discreteness of the spectrum.\medskip\\
\textbf{Proof of Lemma \ref{easy lemma}.} Using that the spectrum
of the 1 dimensional operators:
\begin{equation*}
-\Delta+|x|^\nu, \nu >0,
\end{equation*}
has a strictly positive lowest eigenvalue $\lambda_0(\nu)$, by
scaling (see Eq.(\ref{scale potential})) we get the following
lower bound valid in the sense of quadratic forms on
$H^1(\mathbb{R}^n)$:
\begin{equation*}
-\Delta_{x_1,\cdots,x_n} + |x_1|^{\alpha_1} \prod\limits_{j=2}
^{n} |x_j|^{\alpha_j} \geq -\Delta_{x_2,\dots,x_n}
+\lambda_0(\alpha_1)(\prod\limits_{j=2}^{n}
 |x_j|^{\alpha_j})^{\frac{2}{2+\alpha_1}} .
\end{equation*}
By induction, and symmetrization w.r.t. $x_1,\cdots ,x_n$, it is
easy to show that there exists $n$ positive continuous functions
$f_i:\mathbb{R}\rightarrow \mathbb{R}_+$ with
$\lim\limits_{t\rightarrow \pm \infty} f_i(t)=+\infty$ such that:
\begin{equation*}
-\Delta +\prod\limits_{j=1}^{n} |x_j|^{\alpha_j} \geq
\sum\limits_{j=1}^{n} (-\Delta_{x_j}
+f_j(x_j))=\sum\limits_{j=1}^n T_j.
\end{equation*}
A fortiori, e.g. by a min-max argument\footnote{One could also
exponentiate the functional inequality and take the trace.}, the
spectrum of $H_n^\alpha$ is discrete each of the operators $T_j$
appearing in the r.h.s. being a 1-dimensional Schr\"odinger
operator with confining potential. It is also easy to verify that
$f_j(t)\geq C |t|^{\eta_j}$ for some $\eta_j\in
\mathbb{Q}^{*}_{+}$.  $\hfill{\blacksquare}$
\begin{remark}
\rm{As it was observed in the $2$ dimensional case, the
inequality:
\begin{equation*}
Z_Q(t)= \mathrm{Tr}\, e^{-tH^\alpha_n} \leq \prod\limits_{j=1}^n
\mathrm{Tr}\,e^{-tT_j},
\end{equation*}
can be exploited to get an upper bound, e.g. using Bohr-Sommerfed
quantization conditions for each $T_j$ (see \cite{B-Shu}, chapter
5). But these bounds are not good in view of the results stated in
Theorem 2 and 3.}
\end{remark}

\noindent \textbf{Comments and perspectives:}

\noindent To our knowledge very few things are known about the
interpretation of the results stated in Theorem \ref{main result}
and \ref{main result2} in terms of geometry, physics or dynamical
systems. For example, a geometrical interpretation of the
constants appearing in the short time expansion of
$e^{-tH^{\alpha}_n}$ cannot have a 'classical' meaning (for the
usual symplectic structure of the phase space $T^*\mathbb{R}^n$).
At least these coefficients can be used to construct some (global)
measures on eigenvectors:
\begin{equation*}
w(a) =\lim\limits_{t\rightarrow 0^+} \frac{ \mathrm{Tr} \, A
e^{-tH^{\alpha}_n}} { \mathrm{Tr} \, e^{-tH^{\alpha}_n}},\,
A=\mathrm{Op}(a),\, a\in S^0(\mathbb{R}^{2n}),
\end{equation*}
i.e. the asymptotic results stated in theorems \ref{main result}
and \ref{main result2} provides some normalization factors (in
terms of probability measures) for the statistical distribution of
eigenfunctions in the phase space. Also, because of the strong
singularities of the potential (and of the Liouville measure) on
any hypersurface $\{x_j=0\}$, one could expect some concentration
phenomena for the associated eigenfunctions estimates like it was observed in \cite{BPU,Cam}.\medskip\\
Finally the approach based on the asymptotic behavior of the
quantum partition function $Z_Q(t)$ does not allow to see very
much concerning the classical dynamics generated by the potentials
$V_\alpha$ (when the dynamics is globally defined in the usual
way, meaning that the Hamiltonian vector field has the Lipschitz
regularity). But, as a matter of honesty, the usual semi-classical
methods and their underlying stationary-phase approximations seem
to be inefficient because of the infinite volume of energy
surfaces. \medskip\\

\noindent We will conclude the introduction with the proof of theorems 6 and 7.
Of course, because of the scaling properties of our operators
similar results are valid for Dirichlet-Laplacian in domains
$\Omega^n_\alpha(a)=\{x\in\mathbb{R}^n: |x_1|^{\alpha_1} \cdots
|x_n|^{\alpha_n} <a \}$ for any positive $a$.\\

\noindent \textbf{Proof of theorem 6 and 7:} The proof of these results is straightforward if
we use a sequence of potentials with strictly increasing exponents
$(\alpha)_j=j.\alpha=(j\alpha_1,\cdots,j\alpha_n)$, $j>0$. The
ratios of theorem \ref{main result} satisfy:
\begin{equation*}
d_n((\alpha)_j)=\frac{j\alpha_1+...+j\alpha_{n-1}+2}{2j\alpha_n}\rightarrow
\frac{\alpha_1+...+\alpha_{n-1}}{2\alpha_n} \text{ as }
j\rightarrow +\infty.
\end{equation*}
On the other side, if $j$ tends to infinity by homogenity we
have:
\begin{equation*}
V_{(\alpha)_j} =(V_\alpha)^j \rightarrow \left\{
\begin{matrix}
0 \text{ if } V_\alpha <1,\\
1 \text{ if } V_\alpha =1,\\
+\infty \text{ if } V_\alpha >1.
\end{matrix}
\right.
\end{equation*}
By taking the exponential (the potentials are everywhere positive)
we get the desired result since $e^{-t(-\Delta+ V_{(\alpha)_j})}$
converges to $e^{-t (-\Delta^D_{n,\alpha})}$ as $j\rightarrow \infty$ in the
trace norm. $\hfill{\blacksquare}$\\

\noindent The rest of the paper is organized as follows: in
section 2 we recall the sliced bread estimate used to get good
upper bound on $Z_Q(t)$. Section 3 then contains the proof of
Theorem 2, while in section 4 we prove the technically more
involved asymptotics of Theorem 3.
\section{Slicing techniques for the partition function:}
As was noted in the last section, the Tauberian theorem of Karamata allows us to focus on the small
time divergence of the partition function $Z_Q(t)$. When analyzing
the trace of an operator of this type, it is useful to 'slice' the
problem. That is, we write an operator $A=-\Delta+V(x,y)$, $(x,y) \in \mathbb{R}^a \times
\mathbb{R}^b = \mathbb{R}^c$  on $L^2(\mathbb{R}^c)$, with $V(x,y)$ continuous
and bounded from below as a sum:
\begin{equation*}
A=-\Delta_x + A_x.
\end{equation*}
Here $A_x=-\Delta_y+V(x,y)$ as an operator on $L^2(\mathbb{R}^b)$ depending on $x$.
Let $\lambda_k(x)$ be the increasing sequence of eigenvalues of $A_x$, repeated according to
their multiplicities. Define:
\begin{gather*}
Z_{SB}(t)=\sum_k \mathrm{Tr}_{L^2(\mathbb{R}^{a})}(e^{-t(\Delta_{x}+\lambda_k(x))}),\\
Z_{SGT}(t)= \int \frac{e^{-t||\xi||^2}}{(2\pi)^{a}} \mathrm{Tr}_{L^2(\mathbb{R}^b)}(e^{-t(A_x)})d^{a}\xi d^{a}x,\\
Z_{cl}(t)= \int \frac{e^{-t(||\xi||^2+V(x,y))}}{(2\pi)^{c}} d^{c}\xi d^{a}x d^{b}y.
\end{gather*}
Here SB stands for sliced bread, SGT is sliced Golden-Thompson. The 'sliced-bread' and
'sliced-Golden-Thompson' techniques are now centered around the
following theorem, which is due to B. Simon \cite{Sim2}:
\begin{theorem} \textbf{(Barry Simon's sliced bread inequalities.)}\\
For each $t>0$ we have:
\begin{equation*}
Z_Q(t)\leq Z_{SB}(t) \leq Z_{SGT}(t) \leq Z_{cl}(t).
\end{equation*}
\end{theorem}
The potentials we deal with in this paper provide examples
where $Z_{cl}(t)=\infty$ and, depending on the choice of $\alpha_i$, even $Z_{SGT}(t)=\infty$,
yet the traces $Z_{SB}(t)$ and $Z_Q(t)$ exist.
This provides a set of examples where these estimates prove to be more powerful than the classical one. In the case of the potentials covered in this paper, working with either $Z_{SGT}(t)$ or $Z_{SB}(t)$ leads to studying 'partial trace' functions of the type:
\begin{equation*}
F(x_n,t)=\mathrm{Tr}_{L^2(dx_1...dx_{n-1})}(\exp[-t(\Delta_{x_1...x_{n-1}}+\prod\limits_{i=1}^n
|x_i|^{\alpha_i})])
\end{equation*}
This function satisfies a remarkable functional equation:
\begin{lemma}\label{lemma scaling}
The function $F(x_n,t)$ defined above satisfies the
following scaling relation:
\begin{equation*}
F(x_n,t)=F(x_n t^{d_n},1)=F(1,t |x_n|^{1/d_n}).
\end{equation*}
\end{lemma}
\noindent\textbf{Proof.} The homogeneity of the potential is
crucial in proving this. Use the scaling relations of
Eqs.(\ref{scale potential},\ref{scale laplace}) and thus show:
\begin{gather*}
F(x_n,t)=\mathrm{Tr}_{L^2(dx_1...dx_{n-1})}(\exp[-t(-\Delta_{x_1...x_{n-1}}+\prod\limits_{i=1}^n |x_i|^{\alpha_i})])\\
=\mathrm{Tr}_{L^2(dx_1...dx_{n-1})}(\exp[-(-\Delta_{x_1...x_{n-1}}+ |t^{d_n} x_n|^{\alpha_n} \prod\limits_{i=1}^{n-1} |x_i|^{\alpha_i})])\\
=\mathrm{Tr}_{L^2(dx_1...dx_{n-1})}(\exp[-t  |x_n|^{(1/{d_n})}
(-\Delta_{x_1...x_{n-1}}+  \prod\limits_{i=1}^{n-1}
|x_i|^{\alpha_i})]).
\end{gather*}
For the first equality set $c=t$ and apply (1) then (2), for the second one set $c=|x_n|$ and work in reverse order.
\hfill $\square$
\section{Eigenvalue asymptotics, $\alpha_i \neq \alpha_j, \quad\forall i \neq j$}
This case is techniqually less involved due to the finiteness of $Z_{SGT}(t)$. When all indices
are different, up to permutation of coordinates we can assume
that:
\begin{equation*}
V_\alpha(x)=\prod\limits_{i=1}^n |x_i|^{\alpha_i}, \text{with }
\alpha_1 > ... > \alpha_n.
\end{equation*}
The proof will be by induction over dimension. The dimension 2
case was shown by B. Simon in \cite{Sim2}. Suppose now that in
dimension $n-1$, we have:
\begin{equation*}
\lim_{t\rightarrow 0} t^{(d_{n-1}+1/2)}
\mathrm{Tr}(e^{-tH^{\alpha}_{n-1}})=\mathrm{Tr}((H^{\alpha}_{n-2})^{-d_{n-1}})\pi^{-1/2}
\Gamma(d_{n-1}+1).
\end{equation*}
We will estimate $Z_Q(t)$ from above and below, and then both show
that the difference of these bounds asymptotically goes to zero
and compute the asymptotics of the upper bound. The lower bound
will be found using the Feynman-Kac formula that gives a
representation of the trace $Z_Q(t)$ of the heat kernel as an
expectation value of Brownian motion running for time $2t$. The
upper bound will be found using the sliced Golden-Thompson trace.
We will slice in direction of the coordinate of smallest power in
the potential, $x_n$.
\begin{remark}
\rm{The slicing for $Z_{SGT}(t)$ works
only if one takes slices in the right direction, that is, the one
with smallest exponent in the potential. As we will see one has no choice, as the
integrals will not converge if one slices in a different way.}
\end{remark}
To start, we note that from the sliced bread inequalities, we know already that: $Z_Q(t)\leq Z_{SGT}(t)$. Let us first rewrite this upper bound a little.
Doing the $\xi$-integral explicitly, we get:
\begin{equation*}
Z_{SGT}(t)=(\pi t)^{-1/2} \int_0^{\infty}F(x_n,t)dx_n.
\end{equation*}\\
\noindent \textbf{Lower bound:} Next we will prove a lower bound
that is easy to compare to the expression of $Z_{SGT}$ as an
integral of $F$. Using the Feynman-Kac formula we rewrite $Z_Q(t)$
just like Barry Simon did in the two dimensional case, see
\cite{Sim,Sim2}:
\begin{equation*}
Z_Q(t)= (4\pi t)^{-n/2} \int\limits_{x\in\mathbb{R}^n}
\mathbb{E}_{x,x;2t}[\exp(-\int_0^{2t}\frac{1}{2}
|b_1(s)|^{\alpha_1} ... |b_{n-1}(s)|^{\alpha_{n-1}}
|b_n(s)|^{\alpha_n} ds)]dx,
\end{equation*}
where $b_i$ denotes the 1-dimensional Brownian motion and
$\mathbb{E}_{x,x;2t}$ is the conditional expectation value w.r.t.
the Brownian motion with conditions to start and end\footnote{Of
course, the fact that $\mathbb{E}_{x,y;2t}$ is restricted to the
diagonal $\{x=y\}$ corresponds to the fact that $Z_Q$ is the trace
of the heat kernel.} in $x$ in time $2t$. We proceed by cutting
off all paths such that:
\begin{equation*}
\sup_{0\leq s\leq 2t}|b_n(s)-x_n|>1,
\end{equation*}
and replacing $|b_n(s)|^{\alpha_n}$ by its upper bound
$(|x_n|+1)^{\alpha_n}$. Using that our potential is a product
function, and since the probability measure of the $n$-dimensional
Brownian motion is a product measure, this gives a lower bound for
$Z_Q(t)$, namely:
\begin{gather*}
Z_Q(t)\geq  (\pi t)^{-\frac{1}{2}}(1-\rho (t))\int_0^{\infty}F(|x_n|+1,t)dx_n\\
=(\pi t)^{-\frac{1}{2}}(1-\rho (t))\int_1^{\infty}F(x_n,t)dx_n.
\end{gather*}
This inequality is valid by symmetry and by the following two facts. First,
the probability that a path leaves a compact interval $[x_n-1,x_n+1]$ during a
small time interval $[0,2t]$ is small:
\begin{equation*}
\rho (t) \geq \mathrm{Prob} \sup\limits_{0 \leq s \leq 2t}
(|b_n(s)-x_n|>1),
\end{equation*}
with $\rho (t) \rightarrow 0$ as $t \rightarrow 0^+ $. In fact, a classical result concerning the
$1$-dimensional Brownian motion is that for each $\varepsilon>0$,
there exists a positive constant $C(\varepsilon)$ such that:
\begin{equation*}
\rho (t) \leq C(\varepsilon) e^{-(1-\varepsilon)/4t},\, \text{as }
t\rightarrow 0^+.
\end{equation*}
Second, the monotony of the exponential, the
positivity of the integral and the product structure of the
potential gives a lower bound at $|x_n|+1$.\medskip\\
\noindent \textbf{Upper bound:} Let us now analyze the upper bound
more closely. Theorem 2 follows from the following three
statements:
\begin{gather}
\lim_{t\rightarrow 0} t^{d_n} \int_0^1 F(x_n,t) dx_n = 0 .
\label{Point2}\\
\lim_{t\rightarrow 0} t^{d_n} \int_0^{\infty} F(x_n,t) dx_n =
\mathrm{Tr}((H^\alpha_{n-1})^{-d_n}) \Gamma (d_n+1). \label{Point1}\\
\mathrm{Tr}((H^\alpha_{n-1})^{-d_n})<\infty. \label{Point3}
\end{gather}
This will finish the proof because equation (\ref{Point2}) implies that our upper bound and lower bound agree asymptotically, while equations (\ref{Point1},\ref{Point3}) prove the asymptotic power and constant.
We treat now the result stated in equation (\ref{Point2}). According to
the hypothesis of induction we have:
\begin{equation*}
F(1,t) \sim C t^{-\frac{1}{2}-d_{n-1}},\text{ as } t\rightarrow
0^+,
\end{equation*}
for some $C>0$. It follows that for $|x|$ in a compact subset of $\mathbb{R}_+$
and $t>0$ small:
\begin{equation*}
F(x,t) \leq \tilde{C} (|x|^{\frac{1}{d_n}})^{-\frac{1}{2}-d_{n-1}}
t^{-\frac{1}{2}-d_{n-1}}.
\end{equation*}
Hence since:
\begin{gather*}
d_n -\frac{1}{2}-d_{n-1}=\frac{\alpha_1+\dots
+\alpha_{n-1}+2}{2\alpha_n} -\frac{1}{2} -\frac{\alpha_1+\dots
+\alpha_{n-2}+2}{2\alpha_{n-1}}\\
=\frac{\alpha_1+\dots +\alpha_{n-1}+2}{2\alpha_n} -
\frac{\alpha_1+\dots +\alpha_{n-1}+2}{2\alpha_{n-1}}>0.
\end{gather*}
The decay w.r.t. $t$ is achieved. It remains to show that the
integral w.r.t. $x$ is finite. For $t$ small we have that:
\begin{equation*}
\int\limits_{0}^{1} F(x,t)dx \leq C(t) \int\limits_{0}^{1}
|x|^{\frac{1}{d_n}})^{-\frac{1}{2}-d_{n-1}} dx.
\end{equation*}
We see that this singularity in $x=0$ is integrable if and only if:
\begin{equation*}
-\frac{1}{d_n}({\frac{1}{2}+d_{n-1}})>-1,
\end{equation*}
since we have:
\begin{equation*}
\frac{1}{d_n}({\frac{1}{2}+d_{n-1}})=\frac{2\alpha_n}{\alpha_1+\dots
+\alpha_{n-1}+2} \frac{\alpha_1+\dots
+\alpha_{n-1}+2}{2\alpha_{n-1}}=\frac{\alpha_n}{\alpha_{n-1}} <1,
\end{equation*}
we obtain that the integral is finite.
\begin{remark}\label{Rem critical exponents}
\rm{This result concerning the singularity in $x=0$ is
asymptotically exact as $t\rightarrow 0^+$ so that the condition
$\alpha_n<\alpha_{n-1}$ is necessary to get Eq.(\ref{Point2}).
Hence, like in dimension 2, there is a good way to slice imposed
by the exponents of the potential. One slices always in direction
of the smallest exponent. Also sliced Gordon Thompson is
not working with $\alpha_n=\alpha_{n-1}$ ($\log$ singularity at
the origin) and sliced bread is required.}
\end{remark}
Let us now show the result stated in equation (\ref{Point1}). First,
via the functional equation for $F$, we can scale $t$ out as:
\begin{equation*}
\int_0^{\infty}F(x_n,t)dx_n=\int\limits_0^{\infty} F(x_n
t^{d_n},1)dx_n=t^{-d_n}  \int\limits_0^{\infty}F(x_n,1)dx_n.
\end{equation*}
Using again lemma \ref{lemma scaling}, we can scale out the $n-1$
dimensional operator $H^\alpha_{n-1}$ via:
\begin{gather*}
\int_0^{\infty}F(x_n,1)dx_n =\int\limits_0^{\infty}F(1,x_n^{1/d_n})dx_n\\
=d_n\int\limits_0^{\infty} s^{l_n-1} F(1,s)ds=
d_n \mathrm{Tr}\left(\int_0^{\infty}s^{d_n-1} e^{-sH^\alpha_{n-1}}ds \right)\\
=d_n \Gamma (d_n) \mathrm{Tr}((H^{\alpha}_{n-1})^{-d_n}) =\Gamma(d_n+1)
\mathrm{Tr}((H^{\alpha}_{n-1})^{-d_n}) ,
\end{gather*}
as required.\\

\noindent Finally, to achieve the proof of Theorem \ref{main result} it remains to establish the
finiteness of $\mathrm{Tr}((H^\alpha_{n-1})^{-d_n})$, stated in equation (\ref{Point3}). As shown above, we have:
\begin{gather*}
\mathrm{Tr}((H^\alpha_{n-1})^{-d_n})=\frac{1}{\Gamma(d_n+1)}
\int\limits_0^{\infty}F(1,x_n^{1/d_n})dx_n.
\end{gather*}
We now show that the right hand side is finite. As we just have
seen, for small $x_n$, the induction hypothesis tells us that:
\begin{equation*}
F(1,x_n^{1/d_n}) \sim C
x_n^{\frac{1}{d_n}(-\frac{1}{2}-d_{n-1})}=C
x_n^{-\frac{a_n}{a_{n-1}}},\text{ as }x_n \rightarrow 0^+.
\end{equation*}
Since $\frac{a_n}{a_{n-1}}<1$, the singularity at zero is
integrable. On the other hand, for $x_n$ large we see from the
induction hypothesis that $F(1,x_n^{1/d_n})$ decays as
$e^{-cx_n^{1/{d_n}}}$ and it follows that
$\mathrm{Tr}\,((H^\alpha_{n-1})^{-d_n})<\infty$. $\hfill{\blacksquare}$
\section{Eigenvalue asymptotics, $\alpha_1 =\cdots =\alpha_n=\alpha_0$.}
Consider now the singular case of equal exponents in all
directions, that is:
\begin{equation*}
H^\alpha_n=-\Delta +\prod\limits_{i=1}^n |x_i|^{\alpha_0}, \quad
(x_1,...,x_n) \in \mathbb{R}^n, \quad \alpha_0 \in \mathbb{R}_+^*.
\end{equation*}
We will prove the asymptotics of
$Z_Q(t)=\mathrm{Tr}(e^{tH^\alpha_n})$ again based on an induction
in dimension argument. The structure of the proof is once more to
find an upper and lower bound for $Z_Q(t)$ that asymptotically
agree to first order. We have seen in the last section that
$Z_{SGT}(t)=\infty$ in this situation and therefore, we will need
to exploit the final inequality $Z_Q(t)\leq Z_{SB}(t)$ to get a
finite upper bound. We will employ the Feynman-Kac formula to find
a lower bound for $Z_Q(t)$. To begin with, we recall the claimed
asymptotics:
\begin{equation*}
\lim_{t \rightarrow 0^+} t^{(n/2)+\alpha_0^{-1}}\ln(t^{-1})^{n-1}Z_Q(t)(E)=\frac{\Gamma(1+\alpha_0^{-1})(n/2+\alpha_0^{-1})^{(n-1)}}{\pi^{n/2}(n-1)!}.
\end{equation*}
According to the Tauberian theorem of Karamata this imply:
\begin{equation*}
\lim_{E \rightarrow \infty} E^{-(n/2)+\alpha_0^{-1}}\ln(E)^{-n-1}N_{H^\alpha_n}(E)= \frac{\Gamma(1+\alpha_0^{-1})(n/2+\alpha_0^{-1})^{(n-1)}}{\Gamma(n/2+\alpha_0^{-1})\pi^{n/2}(n-1)!}=:a_n.
\end{equation*}
\medskip

\noindent We proceed now with the computation of the asymptotics
of the aforementioned bounds, starting with the upper
one:\medskip\\
\noindent \textbf{Upper bound.} From the sliced bread
inequalities, we know $Z_Q(t)\leq Z_{SB}(t)$. We compute now the
asymptotics of $Z_{SB}$ directly, making heavy use of scaling
arguments. To start, we prove some auxiliary results in the
following lemma:
\begin{lemma}\label{lemma tech}
Let $A_g=-\Delta_{x_n}+g|x_n|^{\gamma}$, $\gamma>0$. Let
$F^{(\gamma)}_g(t)=\mathrm{Tr}(\mathrm{exp}(-tA_g))$ and $N_{A_g}$
be the associated counting function. Also, denote by
$N_{r}:=N_{B_r}$ the eigenvalue counting function of
$B_r:=-\Delta_{x_1,...,x_{n-1}}+r\prod\limits_i^{n-1}|x_i|^{\alpha_0}$.
Omit the indices $g,r$ whenever $r=1, g=1$. Then:
\begin{itemize}
\item[a)] $F^{(\gamma)}_g(t)=F^{(\gamma)}(g^{\tau}t)$, where $\tau=2/(\gamma + 2)$,
\item[b)] $N_r (E) = N(r^{\eta}E)$, where $\eta=((n-1)\alpha_0+2)/2$,
\item[c)]$\lim_{t \rightarrow 0^+} t^{\mu}F^{(\gamma)}(t)=\pi^{-1/2} \Gamma(1+\gamma^{-1})$, where $\mu=(\gamma+2)/2\gamma$,
\item[d)]$\lim_{t \rightarrow 0^+} t^{\mu+1}F^{\prime(\gamma)}(t)=-\mu \pi^{-1/2} \Gamma(1+\gamma^{-1})$.
\end{itemize}
\end{lemma}
\noindent \textbf{Proof:} Part a) and b) follow by the scaling
relations (1) and (2). Part c) follows from the observation that
since for the operator $A_g$ in one dimension the potential
$|x_n|^\gamma$ is confining, we actually have
$\frac{Z_{cl}(t)}{F^{(\gamma)}(t)} \to 1$ as $t \to 0^+$, cf.
\cite{Sim}. Now:
\[Z_{cl}(t)=\frac{1}{2\pi}\int\limits_{\mathbb{R}^2} e^{-t(\xi^2+|x|^\gamma)}dxd\xi=t^{-\gamma}\pi^{-1/2}\Gamma(1+\gamma^{-1}),\]
and the claim follows. Finally, part d) which tells us that the
asymptotics of the differential of $F$ are given as the
differential of the leading term of $F$ is shown by the following
calculation:
\begin{eqnarray*}
-F^{\prime(\gamma)}(t)&=&\int\limits_0^\infty xe^{-tx}dN_{A}(x)\\
&=&\int\limits_0^\infty e^{-tx}(xt-1)N_A(x)dx\\
&=&t^{-1}\int\limits_0^\infty e^{-y}(y-1)N_A(\frac{y}{t})dy.
\end{eqnarray*}
To conclude that the claimed asymptotics of $F^{\prime(\gamma)}$ hold we use the Karamata Tauberian theorem for $N_A$ and apply result $c)$.
 $\hfill \square$.\\

\noindent With the help of this lemma, we first rescale the sliced bread trace $Z_{SB}(t)$ and then represent
it as an integral expression which is suitable for the computation of the asymptotics.
Since $N_r$ is, up to the factor $r$ in the potential, the counting function of the operator
$B_1=H_{n-1}^\alpha$ in dimension $n-1$, this will get us in position to use the induction hypothesis.\\

\noindent For the rescaling, let $\epsilon_j(x_n)$ be the $j$th eigenvalue of $-\Delta_{x_1,...,x_{n-1}}+|x_n|^{\alpha_0} \prod\limits_{i=1}^{n-1}|x_i|^{\alpha_0}$ so with $d_n=\frac{n-1}{2}+\alpha_0^{-1}$ we get by relation b):

\[ \epsilon_j(x_n)=|x_n|^{1/{d_n}}\epsilon_j(1)=:|x|^{1/{d_n}}\epsilon_j.\]
Now compute:
\begin{eqnarray*}
Z_{SB}(t)&=&\sum\limits_j \mathrm{Tr} (\exp [-t(\Delta_x + \epsilon_j(x_n))])\\
&=&\sum\limits_j \mathrm{Tr} (\exp [-t(\Delta_x + \epsilon_j|x_n|^{1/d_n})])\\
&=&\sum\limits_j \mathrm{Tr} (\exp [-t\epsilon_j^{b_n}(\Delta_x + |x_n|^{1/d_n})]).
\end{eqnarray*}
where we have used scaling relation b) so:  $b_n=\frac{d_n}{d_n+(1/2)}$. We now represent $Z_{SB}$ as an integral and use integration by parts:
\begin{eqnarray*}
Z_{SB}(t)&=&\sum\limits_j F^{(1/d_n)}(t\epsilon_j^{b_n})\\
&=&\int\limits_0^\infty  F^{(1/d_n)}(tE^{b_n})dN(E) \\
&=&-\int\limits_0^{\infty}tbE^{b_n-1} F^{\prime(1/d_n)}(tE^{b_n})N(E)dE.
\end{eqnarray*}
There are no boundary terms for the following reasons:\\
First, we have $N=0$ for small $E>0$ since a direct evaluation shows that $0$ is no eigenvalue of $H^\alpha_{n-1}$ and the spectrum of this operator is discrete by the hypothesis of induction.\\
Second, $F^{(1/d_n)}N \rightarrow 0$ for fixed $t$ and $E\rightarrow \infty$, since for large energies, the trace $F^{(1/d_n)}$ decays exponentially.\\

\noindent The rescaling led to $tE^{b_n}$ as the argument of $F^{\prime(1/d_n)}$. This will prove to be crucial for the precise asymptotic analysis. We proceed from here very similarly to Simon's arguments for dimension 2 by analyzing this integral on several pieces. As we just argued, since $N(E)=0$ for $E>0$ small by our hypothesis, it is sufficient to integrate from the ground state energy $E_0>0$ to $\infty$. Now we pick values $E_0<E_1<E_2<\infty$ such that:

\begin{equation*}
E_1^{b_n}t=|\ln t|^{-1}, \quad \quad E_2^{b_n}t=1.
\end{equation*}
We will now estimate the integral on $(E_0,E_1)$, $(E_1,E_2)$, $(E_2,\infty)$ separately and see that only the integral
on $(E_0,E_1)$ will contribute to the top order coefficients.\\

\begin{itemize}
\item[$(E_2,\infty)$:]
On this piece we first note that $N(E)\leq c E^{d_n}\ln(E)^{n-2}$ for all $E$ because of the hypothesis of induction, and in addition, in the region $y\geq 1$, $-F^{\prime(1/d_n)}(y)=\sum \tilde{\epsilon}_je^{-y\tilde{\epsilon}_j}\leq D e^{-cy}$, with $D,C>0$ and $ \tilde{\epsilon}_j$ being the eigenvalues of $\Delta_{x_n}+|y|^{1/d_n}$. Thus we estimate:

\begin{eqnarray*}
&&-\int_{E_2}^{\infty} (tb_n)E^{b_n-1}N(E)F^{\prime(1/d_n)}(tE^{b_n})dE\\ &\leq& c_1 \int_{E_2}^{\infty}t E^{b_n-1+d_n}\ln(E)^{n-2}\exp(-ctE^{b_n})dE\\
&=&c_2t^{-b_n/d_n}\int_1^{\infty} (\ln\frac{y}{t})^{n-2} y^{d_n/b_n}e^{-cy}dy,
\end{eqnarray*}
for some constants $c_1$, $c_2$. We see that all terms in the integral are bounded by $t^{-d_n/b_n}\ln(t^{-1})^{-(n-2)}=t^{-(d_n+1/2)}\ln(t^{-1})^{-(n-2)}$ which is small on the level of $t^{-(d_n+1/2)}\ln(t^{-1})^{-(n-1)}$.\\

\item[$(E_1,E_2)$:] We bound $N(E)$ as before but use now
$-F^{\prime(1/d_n)}(y)\leq Cy^{-d_n -\frac{3}{2}}$, by result d)
of the previous lemma. Then, where c is a constant which changes
from equation to equation:
\begin{eqnarray*}
&&-\int_{E_1}^{E_2} (tb_n)E^{b_n-1}N(E)F^{\prime(1/d_n)}(tE^{b_n})dE \\ &\leq& c \int_{E_1}^{E_2}t E^{b_n-1+d_n}\ln(E)^{n-2}(tE^{b_n})^{-d_n-3/2}dE\\
&\leq&ct^{-d_n-1/2} \int_{E_1}^{E_2}\ln(E)^{n-2}E^{-1}dE\\
&\leq& ct^{-d_n-1/2} \ln(E_2/E_2)^{n-1},
\end{eqnarray*}
where we have used $b_n-1+d_n-b_n(d_n+\frac{3}{2})=-1$, since $b_n(d_n+\frac{1}{2})=d_n$. Because $\ln(E_2/E_1)=b_n^{-1}\ln(|\ln t|)$, this integral is small compared to $t^{-(d_n+1/2)}\ln(t^{-1})^{-(n-1)}$.\\

\item[$(E_0,E_1)$:]
Finally, for the last piece we replace $F^{\prime(1/d_n)}$ by its asymptotic value making a multiplicative error of the form $1+o(1)$. In other words, we can bound $F'$ from above and below by:\\
\[-(1\pm \epsilon(t))(d_n+\frac{1}{2}) \pi^{-1/2} \Gamma(1+d_n)(tE^{b_n})^{(d_n+3/2)},\]
with $\epsilon(t) \rightarrow 0$. This is true because the argument of $F'$ namely $tE^{b_n}$ is now bounded from above by $|\ln(t)|^{-1}$, which implies:
\[ |F^{\prime(1/d_n)}(tE^{b_n})-(d_n+\frac{1}{2}) \pi^{-1/2} \Gamma(1+d_n)(tE^{b_n})^{(d_n+3/2)}|=o(\frac{1}{|\ln(t)|}).\]
Thus, if $\sim$ means the ratio goes to 1, we see that:
\begin{eqnarray*}
A(t)&:=&-\int_{E_0}^{E_1} (tb_n)E^{b_n-1}N(E)F^{\prime(1/d_n')}(tE^{b_n})dE \\
&\sim&\int_{E_0}^{E_1}tb_n E^{-1}\ln(E)^{n-2}t^{-d_n-3/2}(d_n+1/2) \pi^{-1/2}\Gamma(1+d_n) a_{n-1}\\
&&\left[ \frac{N(E)}{a_{n-1}E^{d_n}\ln(E)^{n-2}}\right]dE,
\end{eqnarray*}
where $a_{n-1}$ is just the constant of the eigenvalue counting function in dimension $(n-1)$. Now observe that as $t \to 0^+$ the value of $E_1$ tends to $\infty$, while $E_0$ remains constant. Since by the induction hypothesis we also have:
\[\frac{N(E)}{a_{n-1}E^{d_n}\ln(E)^{n-2}} \rightarrow 1 \quad \mathrm{as}~ E \to \infty
,\]
we see now that we have another asymptotic equivalence, namely:
\begin{eqnarray*}
A(t)&\sim& a_{n-1}b_n (d_n+1/2) \pi^{-1/2}\Gamma(1+d_n)  t^{-d_n-1/2} \int_{E_0}^{E_1} E^{-1}\ln(E)^{n-2}dE\\
&=&\frac{a_{n-1}d_n\Gamma(1+d_n)}{\pi^{1/2}}  t^{-d_n-1/2}\frac{\ln(E_1/E_0)^{n-1}}{(n-1)},
\end{eqnarray*}
Now we compute:
\begin{equation*}
\ln(E_1/E_0)=\ln[ct^{-1/b_n}|\ln t|^{-1/b_n}] \sim \frac{1}{b_n}\ln(t^{-1}).
\end{equation*}
where $c$ is once more some constant. We conclude that:
\begin{eqnarray*}
A(t)t^{-d_n-1/2}\ln(t^{-1})^{-(n-1)}&\sim& \pi^{-1/2}a_{n-1}d_n b_n^{-(n-1)}\Gamma(1+d_n)(n-1)^{-1}\\
&=&\frac{\Gamma(d_3) (d_{n+1})^{n-1}}{\pi^{n/2}(n-1)!}.
\end{eqnarray*}
So we get the claimed asymptotics for the upper bound:
\begin{eqnarray*}
\lim\limits_{t\rightarrow 0} t^{n/2+\alpha_0^{-1}}\ln(t^{-1})^{-(n-1)}\mathrm{Tr}(e^{-tH^\alpha_n})\leq\frac{\Gamma(1+\alpha_0^{-1}) (n/2+\alpha_0^{-1})^{n-1}}{\pi^{n/2}(n-1)!}.
\end{eqnarray*}
\end{itemize}
\medskip

\noindent \textbf{Lower bound.}
For the lower bound we will employ the Feynman-Kac formula again. Before we do so we prove the following lemma:
\begin{lemma}
Let $f: \mathbb{R}\rightarrow \mathbb{R}$ be any integrable function. Let $P:\mathbb{R}^n \to \mathbb{R}$ denote the map given by $P(x)=\prod\limits_{i=1}^n|x_i|$. Let $X_{a}:=\{x \in \mathbb{R}^n~|~ \exists x_i: |x_i|<a \}$. Then we have:

\begin{equation*}
\int_{\mathbb{R}^n \setminus X_{a}} f\circ P ~dx = 2^n \int_{a^n }^\infty \frac{f(p)}{(n-1)!}\log(\frac{p}{a^n})^{n-1} dp.
\end{equation*}

\end{lemma}

\paragraph{Proof:}
We prove this by induction. If n=1 the equality is easily verified. Now:
\begin{eqnarray*}
&&\int_{\mathbb{R}^n \setminus X_{a}} f \circ P ~ dx\\
&=& 2^n \int_{a}^\infty \left(\int_{a}^\infty \cdots \int_{a}^\infty f(x_1\cdot ... \cdot x_n) \right) dx_1 \cdots dx_{n-1}) dx_n\\
&=& 2^n \int_{a}^\infty \left(\int_{a^{n-1}}^\infty \frac{f(x_n\cdot q)}{(n-2)!} \log(\frac{q}{a^{n-1}})^{n-2} dq \right) dx_n\\
&=& 2^n \int_{a^n}^\infty f(p)\left(\int_{a}^{\frac{p}{a^{n-1}}} \frac{1}{(n-2)!} \left[\log(\frac{p}{x_na^{n-1}})\right]^{n-2} \frac{1}{x_n} dx_n\right) dp\\
&=& 2^n \int_{a^n}^\infty f(p)\left( \frac{1}{(n-2)!} \left[\sum_{k=0}^{n-2}\binom{n-2}{k}\log(\frac{p}{a^{n-1}})^k\int_{a}^{\frac{p}{a^{n-1}}}(-\log(x_n))^{n-k-2}\frac{1}{x_n}dx_n\right] \right) dp\\
&=& 2^n \int_{a^n}^\infty \frac{f(p)}{(n-1)!} \left(  \left[ - \sum_{k=0}^{n-2}\binom{n-1}{k} \log(\frac{p}{a^{n-1}})^k (-\log(\frac{p}{a^{n-1}}))^{n-k-1}\right] \right. \\
&& + \left. \left[ \sum_{k=0}^{n-2}\binom{n-1}{k}\log(\frac{p}{a^{n-1}})^k (-\log(a))^{n-k-1}\right] \right) dp\\
&=& 2^n \int_{a^n}^\infty \frac{f(p)}{(n-1)!} \left( - \left[ \log(\frac{p}{a^{n-1}}) - \log(\frac{p}{a^{n-1}})^{n-1} +  \log(\frac{p}{a^{n-1}})^{n-1}\right] \right. \\
&& + \left. \left[ \log(\frac{p}{a^{n-1}}) - \log(a)^{n-1} - \log(\frac{p}{a^{n-1}})^{n-1} \right]\right)dp\\
&=& 2^n \int_{a^n}^\infty \frac{f(p)}{(n-1)!}\log(\frac{p}{a^n})^{n-1} dp.
\end{eqnarray*}$\hfill \square$\\

\noindent We compute the lower bound as follows. Apply the Feynman-Kac formula for $Z_Q(t)$ to obtain:
\begin{equation*}
Z_Q(t)= (4\pi t)^{-n/2} \int\limits_{x\in\mathbb{R}^n}
\mathbb{E}_{x,x;2t}[\exp(-\int_0^{2t}\frac{1}{2}\prod\limits_{i=1}^n
|b_i(s)|^{\alpha_0}ds)]dx,
\end{equation*}
Now, remove some of the points $x=(x_1,...,x_n)$ in the integral, as well as some Brownian paths. Only keep those points $x$ that satisfy $x \notin X_{t^{1/2}(\ln t)^2}$. Only keep Brownian paths with $\sup_{0\leq s\leq 2t} |b_i(s)-x_i|\leq t^{1/2}|\ln(t)|$, $\forall 1\leq i\leq n$, where we remind again that in our notation $b_i(0)=x_i$.\\

\noindent The measure of such paths is $1-\rho(t)$, where $\rho(t)\sim e^{-D(\ln(t)^2)}$ as $t \to 0^+$. Since $P(x)=\prod\limits_{i=1}^n |x_i|$, we have on points $x \in \mathbb{R}^n$ with $x_i \neq 0$, $\forall i$:
$$\frac{\partial \ln P(x)}{\partial x_i}= \frac{1}{x_i}, \quad \forall i.$$
Therefore, everywhere along the remaining paths we get by Taylor's formula:
\begin{eqnarray*}|\ln P(b(t))-\ln P(x)|&=&|D(\ln P)_{|x}(b(t)-x)+ o(||b(t)-x||)|\\
&\leq& c/|\ln(t)|.
\end{eqnarray*}

\noindent Thus, with $\kappa(t):=\exp(c/\ln(t))$ we have:
$P(b(t))\leq \kappa(t)P(x)$. Inserting this into the expectation
value, the integrand in the $s$-integral is no longer dependant on
the Brownian motion. Thus:
\begin{eqnarray*}
\mathrm{Tr}(e^{-tH^\alpha_n})&\geq&\frac{1-\rho(t)}{(4\pi t)^{n/2}}\int\limits_{\mathbb{R}^n \setminus X_{t^{1/2}\ln(t)^2}} e^{-t\kappa^{\alpha_0}P(x_1,...,x_n)^{\alpha_0}}dx_1\cdots dx_n\\
&=&\frac{1-\rho(t)}{(\pi t)^{n/2}}\int\limits_{t^{n/2}\ln(t)^{2n}}^{\infty} \frac{e^{-t\kappa^{\alpha_0}p^{\alpha_0}}}{(n-1)!}\ln(\frac{p}{t^{n/2}\ln(t)^{2n}})^{n-1}dp\\
&=&\frac{1-\rho(t)}{\kappa t^{\alpha_0^{-1}}(\pi t)^{n/2}}\int\limits_{t^{n/2+\alpha_0^{-1}}\ln(t)^{2n}\kappa}^{\infty} \frac{e^{-w^{\alpha_0}}}{(n-1)!}\ln(\frac{w}{\kappa t^{n/2+\alpha_0^{-1}}\ln(t)^{2n}})^{n-1}dw\\
\end{eqnarray*}
\noindent Since $\rho\rightarrow 0$ and $\kappa \rightarrow 1$ as
$t \to 0^+$ we obtain:
\begin{eqnarray*}
\liminf\limits_{t\rightarrow 0} t^{n/2+\alpha_0^{-1}}(\ln(t^{-n/2-\alpha_0^{-1}}))^{-(n-1)}\mathrm{Tr}(e^{-tH^\alpha_n})\geq\frac{1}{\pi^{n/2}(n-1)!}\int\limits_0^{\infty}e^{-w^{\alpha_0}}dw.
\end{eqnarray*}
\noindent The Integral is equal to
$(\alpha_0^{-1})\Gamma(\alpha_0^{-1})=\Gamma(1+\alpha_0^{-1})$ and
since
$\ln(t^{-n/2-\alpha_0^{-1}})=(n/2+\alpha_0^{-1})\ln(t^{-1})$, we
now conclude:

\begin{eqnarray*}
\liminf\limits_{t\rightarrow 0} t^{n/2+\alpha_0^{-1}}(\ln(t^{-1}))^{-(n-1)}\mathrm{Tr}(e^{-tH^\alpha_n})\geq \frac{\Gamma(1+\alpha_0^{-1})(n/2+\alpha_0^{-1})^{(n-1)}}{\pi^{n/2}(n-1)!}.
\end{eqnarray*}

\noindent So lower and upper bound agree to first order and the
theorem is proven.$\hfill \square$\medskip

\noindent \textbf{Acknowledgments.} The authors thank Werner
Kirsch and Brice Franke for explaining us Feynman-Kac
representations of heat-kernels/operators and the Ito-calculus.
Both authors were supported by the project \textit{SFB-TR12},
\textit{Symmetries and Universality in Mesoscopic Systems} founded
by the DFG.

\end{document}